\documentclass[12pt]{article}

\usepackage{amsmath,amsthm}
\usepackage{amssymb}
\usepackage{color}
\usepackage{xspace}
\usepackage[colorlinks=true,
linkcolor=webgreen,
filecolor=webbrown,
citecolor=webgreen]{hyperref}

\definecolor{webgreen}{rgb}{0,.5,0}

\def\blue{\textcolor{blue} }
\def\v{\vert}

\def\u{\ensuremath{\mathcal U}\xspace}

\def\mbf#1{\mathchoice{\hbox{\boldmath $\displaystyle #1$}}
        {\hbox{\boldmath $\textstyle #1$}}
        {\hbox{\boldmath $\scriptstyle #1$}}
        {\hbox{\boldmath $\scriptscriptstyle #1$}}} 

\newcommand{\seqnum}[1]{\href{http://www.research.att.com/cgi-bin/access.cgi/as/~njas/sequences/eisA.cgi?Anum=#1}{\underline{#1}}}

\catcode`\@=11

\thicklines
\newskip\Einheit \Einheit=.6cm
\newcount\xcoord \newcount\ycoord
\newdimen\xdim \newdimen\ydim \newdimen\PfadD@cke \newdimen\Pfadd@cke
\def\PfadDicke#1{\PfadD@cke#1 \divide\PfadD@cke by2 
\Pfadd@cke\PfadD@cke \multiply\PfadD@cke by2}
\long\def\LOOP#1\REPEAT{\def\BODY{#1}\ITERATE}
\def\ITERATE{\BODY \let\next\ITERATE \else\let\next\relax\fi \next}
\let\REPEAT=\fi
\def\Punkt{\hbox{\raise-2pt\hbox to0pt{\hss\scriptsize$\bullet$\hss}}}

\def\DuennPunkt(#1,#2){\unskip
  \raise#2 \Einheit\hbox to0pt{\hskip#1 \Einheit
          \raise-1.5pt\hbox to0pt{\hss\tiny$\bullet$\hss}\hss}}
		  
\def\NormalPunkt(#1,#2){\unskip
  \raise#2 \Einheit\hbox to0pt{\hskip#1 \Einheit
          \raise-3pt\hbox to0pt{\hss\large$\bullet$\hss}\hss}}
\def\DickPunkt(#1,#2){\unskip
  \raise#2 \Einheit\hbox to0pt{\hskip#1 \Einheit
          \raise-4pt\hbox to0pt{\hss\Large$\bullet$\hss}\hss}}
\def\Kreis(#1,#2){\unskip
  \raise#2 \Einheit\hbox to0pt{\hskip#1 \Einheit
          \raise-4pt\hbox to0pt{\hss\Large$\circ$\hss}\hss}}
\def\Diagonale(#1,#2)#3{\unskip\leavevmode
  \xcoord#1\relax \ycoord#2\relax
      \raise\ycoord \Einheit\hbox to0pt{\hskip\xcoord \Einheit
         \unitlength\Einheit
         \line(1,1){#3}\hss}}
\def\AntiDiagonale(#1,#2)#3{\unskip\leavevmode
  \xcoord#1\relax \ycoord#2\relax \advance\xcoord by -0.05\relax
      \raise\ycoord \Einheit\hbox to0pt{\hskip\xcoord \Einheit
         \unitlength\Einheit
         \line(1,-1){#3}\hss}}
\def\Pfad(#1,#2),#3\endPfad{\unskip\leavevmode
  \xcoord#1 \ycoord#2 \thicklines\ZeichnePfad#3\endPfad\thinlines}
\def\ZeichnePfad#1{\ifx#1\endPfad\let\next\relax
  \else\let\next\ZeichnePfad
    \ifnum#1=1
      \raise\ycoord \Einheit\hbox to0pt{\hskip\xcoord \Einheit
         \vrule height\Pfadd@cke width1 \Einheit depth\Pfadd@cke\hss}%
      \advance\xcoord by 1
     \else\ifnum#1=2
      \raise\ycoord \Einheit\hbox to0pt{\hskip\xcoord \Einheit
         \unitlength\Einheit
         \line(0,1){1}\hss}
      \advance\xcoord by 0
      \advance\ycoord by 1
 \else\ifnum#1=3
      \raise\ycoord \Einheit\hbox to0pt{\hskip\xcoord \Einheit
         \unitlength\Einheit
         \line(1,1){1}\hss}
      \advance\xcoord by 1
      \advance\ycoord by 1
    \else\ifnum#1=4
      \raise\ycoord \Einheit\hbox to0pt{\hskip\xcoord \Einheit
         \unitlength\Einheit
         \line(1,-1){1}\hss}
      \advance\xcoord by 1
      \advance\ycoord by -1
   \else\ifnum#1=5
      \raise\ycoord \Einheit\hbox to0pt{\hskip\xcoord \Einheit
         \unitlength\Einheit
         \line(2,1){2}\hss}
      \advance\xcoord by 2
      \advance\ycoord by 1
	  \else\ifnum#1=6
      \raise\ycoord \Einheit\hbox to0pt{\hskip\xcoord \Einheit
         \unitlength\Einheit
         \line(2,-1){2}\hss}
      \advance\xcoord by 2
      \advance\ycoord by -1
	  \else\ifnum#1=7
      \raise\ycoord \Einheit\hbox to0pt{\hskip\xcoord \Einheit
         \unitlength\Einheit
         \line(3,1){3}\hss}
      \advance\xcoord by 3
      \advance\ycoord by 1
	  \else\ifnum#1=8
      \raise\ycoord \Einheit\hbox to0pt{\hskip\xcoord \Einheit
         \unitlength\Einheit
         \line(3,-1){3}\hss}
      \advance\xcoord by 3
      \advance\ycoord by -1
    \fi\fi\fi\fi\fi\fi\fi\fi
  \fi\next}
\def\hSSchritt{\leavevmode\raise-.4pt\hbox 
to0pt{\hss.\hss}\hskip.2\Einheit
  \raise-.4pt\hbox to0pt{\hss.\hss}\hskip.2\Einheit
  \raise-.4pt\hbox to0pt{\hss.\hss}\hskip.2\Einheit
  \raise-.4pt\hbox to0pt{\hss.\hss}\hskip.2\Einheit
  \raise-.4pt\hbox to0pt{\hss.\hss}\hskip.2\Einheit}
\def\vSSchritt{\vbox{\baselineskip.2\Einheit\lineskiplimit0pt
\hbox{.}\hbox{.}\hbox{.}\hbox{.}\hbox{.}}}
\def\DSSchritt{\leavevmode\raise-.4pt\hbox to0pt{%
  \hbox to0pt{\hss.\hss}\hskip.2\Einheit
  \raise.2\Einheit\hbox to0pt{\hss.\hss}\hskip.2\Einheit
  \raise.4\Einheit\hbox to0pt{\hss.\hss}\hskip.2\Einheit
  \raise.6\Einheit\hbox to0pt{\hss.\hss}\hskip.2\Einheit
  \raise.8\Einheit\hbox to0pt{\hss.\hss}\hss}}
\def\dSSchritt{\leavevmode\raise-.4pt\hbox to0pt{%
  \hbox to0pt{\hss.\hss}\hskip.2\Einheit
  \raise-.2\Einheit\hbox to0pt{\hss.\hss}\hskip.2\Einheit
  \raise-.4\Einheit\hbox to0pt{\hss.\hss}\hskip.2\Einheit
  \raise-.6\Einheit\hbox to0pt{\hss.\hss}\hskip.2\Einheit
  \raise-.8\Einheit\hbox to0pt{\hss.\hss}\hss}}
\def\SPfad(#1,#2),#3\endSPfad{\unskip\leavevmode
  \xcoord#1 \ycoord#2 \ZeichneSPfad#3\endSPfad}
\def\ZeichneSPfad#1{\ifx#1\endSPfad\let\next\relax
  \else\let\next\ZeichneSPfad
    \ifnum#1=1
      \raise\ycoord \Einheit\hbox to0pt{\hskip\xcoord \Einheit
         \hSSchritt\hss}%
      \advance\xcoord by 1
    \else\ifnum#1=2
      \raise\ycoord \Einheit\hbox to0pt{\hskip\xcoord \Einheit
        \hbox{\hskip-2pt \vSSchritt}\hss}%
      \advance\ycoord by 1
    \else\ifnum#1=3
      \raise\ycoord \Einheit\hbox to0pt{\hskip\xcoord \Einheit
         \DSSchritt\hss}
      \advance\xcoord by 1
      \advance\ycoord by 1
    \else\ifnum#1=4
      \raise\ycoord \Einheit\hbox to0pt{\hskip\xcoord \Einheit
         \dSSchritt\hss}
      \advance\xcoord by 1
      \advance\ycoord by -1
    \fi\fi\fi\fi
  \fi\next}
\def\Koordinatenachsen(#1,#2){\unskip
 \hbox to0pt{\hskip-.5pt\vrule height#2 \Einheit width.5pt depth1 
\Einheit}%
 \hbox to0pt{\hskip-1 \Einheit \xcoord#1 \advance\xcoord by1
    \vrule height0.25pt width\xcoord \Einheit depth0.25pt\hss}}
\def\Koordinatenachsen(#1,#2)(#3,#4){\unskip
 \hbox to0pt{\hskip-.5pt \ycoord-#4 \advance\ycoord by1
    \vrule height#2 \Einheit width.5pt depth\ycoord \Einheit}%
 \hbox to0pt{\hskip-1 \Einheit \hskip#3\Einheit 
    \xcoord#1 \advance\xcoord by1 \advance\xcoord by-#3 
    \vrule height0.25pt width\xcoord \Einheit depth0.25pt\hss}}
\def\Gitter(#1,#2){\unskip \xcoord0 \ycoord0 \leavevmode
  \LOOP\ifnum\ycoord<#2
    \loop\ifnum\xcoord<#1
      \raise\ycoord \Einheit\hbox to0pt{\hskip\xcoord 
\Einheit\Punkt\hss}%
      \advance\xcoord by1
    \repeat
    \xcoord0
    \advance\ycoord by1
  \REPEAT}
\def\Gitter(#1,#2)(#3,#4){\unskip \xcoord#3 \ycoord#4 \leavevmode
  \LOOP\ifnum\ycoord<#2
    \loop\ifnum\xcoord<#1
      \raise\ycoord \Einheit\hbox to0pt{\hskip\xcoord 
\Einheit\Punkt\hss}%
      \advance\xcoord by1
    \repeat
    \xcoord#3
    \advance\ycoord by1
  \REPEAT}
\def\Label#1#2(#3,#4){\unskip \xdim#3 \Einheit \ydim#4 \Einheit
  \def\lo{\advance\xdim by-.5 \Einheit \advance\ydim by.5 \Einheit}%
  \def\llo{\advance\xdim by-.25cm \advance\ydim by.5 \Einheit}%
  \def\loo{\advance\xdim by-.5 \Einheit \advance\ydim by.25cm}%
  \def\o{\advance\ydim by.25cm}%
  \def\ro{\advance\xdim by.5 \Einheit \advance\ydim by.5 \Einheit}%
  \def\rro{\advance\xdim by.25cm \advance\ydim by.5 \Einheit}%
  \def\roo{\advance\xdim by.5 \Einheit \advance\ydim by.25cm}%
  \def\l{\advance\xdim by-.30cm}%
  \def\r{\advance\xdim by.30cm}%
  \def\lu{\advance\xdim by-.5 \Einheit \advance\ydim by-.6 \Einheit}%
  \def\llu{\advance\xdim by-.25cm \advance\ydim by-.6 \Einheit}%
  \def\luu{\advance\xdim by-.5 \Einheit \advance\ydim by-.30cm}%
  \def\u{\advance\ydim by-.30cm}%
  \def\ru{\advance\xdim by.5 \Einheit \advance\ydim by-.6 \Einheit}%
  \def\rru{\advance\xdim by.25cm \advance\ydim by-.6 \Einheit}%
  \def\ruu{\advance\xdim by.5 \Einheit \advance\ydim by-.30cm}%
  #1\raise\ydim\hbox to0pt{\hskip\xdim
     \vbox to0pt{\vss\hbox to0pt{\hss$#2$\hss}\vss}\hss}%
}
\catcode`\@=12

\begin{document}

\begin{center}
\vskip 1cm{\LARGE\bf 
Sets, Lists and Noncrossing Partitions     
}
\vskip 1cm
\large
David Callan\\
Department of Statistics \\
University of Wisconsin-Madison  \\
1300 University Ave  \\
Madison, WI \ 53706-1532  \\
USA \\
{\tt callan@stat.wisc.edu}
\end{center}

\vskip .2 in

\begin{abstract}
Partitions of $[n]=\{1,2,\ldots,n\}$ into sets of lists 
(\htmladdnormallink{A000262}{http://www.research.att.com:80/cgi-bin/access.cgi/as/njas/sequences/eisA.cgi?Anum=A000262}) 
are somewhat less numerous than partitions of $[n]$ into lists of sets
(\htmladdnormallink{A000670}{http://www.research.att.com:80/cgi-bin/access.cgi/as/njas/sequences/eisA.cgi?Anum=A000670}).
Here we observe that the former are actually equinumerous with 
partitions of $[n]$ into
lists of \emph{noncrossing} sets and give a bijective proof. We show 
that partitions of $[n]$ into sets of noncrossing lists are counted by 
\htmladdnormallink{A088368}{http://www.research.att.com:80/cgi-bin/access.cgi/as/njas/sequences/eisA.cgi?Anum=A088368} and 
generalize this result to introduce a transform on integer 
sequences that we dub the ``noncrossing partition'' transform.
We also derive recurrence relations to count partitions of $[n]$ into lists of 
noncrossing lists.
\end{abstract}

\vspace{6mm}

\section{Introduction} A \emph{partition} of 
$[n]=\{1,2,\ldots,n\}$ is a collection of nonempty disjoint sets, 
called  \emph{blocks}, whose union is $[n]$. The notion of partition 
can be 
generalized by taking into account the order of the elements within 
each block or the order of the blocks themselves or both. To 
distinguish cases we use the terms list and set with their usual 
connotations of ordered and unordered respectively. Thus there are 
four cases: sets of sets (ordinary set partitions), sets of lists, lists 
of sets, and lists of lists. For unrestricted partitions the four 
counting sequences are respectively the
\htmladdnormallink{Bell numbers}{http://www.research.att.com:80/cgi-bin/access.cgi/as/njas/sequences/eisA.cgi?Anum=A000110},
\htmladdnormallink{A000262}{http://www.research.att.com:80/cgi-bin/access.cgi/as/njas/sequences/eisA.cgi?Anum=A000262},
\htmladdnormallink{A000670}{http://www.research.att.com:80/cgi-bin/access.cgi/as/njas/sequences/eisA.cgi?Anum=A000670},
and
\htmladdnormallink{A002866}{http://www.research.att.com:80/cgi-bin/access.cgi/as/njas/sequences/eisA.cgi?Anum=A002866}.

A partition is \emph{noncrossing} if there do not exist four distinct 
elements $a<b<c<d$ with $a,c$ both in one block and $b,d$ both in 
another. It is well known that noncrossing partitions of $[n]$ (sets 
of noncrossing sets) are counted by the Catalan number $C_{n}$ 
(\htmladdnormallink{A000108}{http://www.research.att.com:80/cgi-bin/access.cgi/as/njas/sequences/eisA.cgi?Anum=A000108}).

In \S 2 we show that partitions of $[n]$ into lists
of noncrossing sets are equinumerous with partitions of $[n]$ into 
arbitrary sets of lists. In \S 3 we show that the ``set of noncrossing 
lists'' case has a generating function 
$A(x)=1+x+3x^{2}+13x^{3}+69x^{4}+\cdots $ that satisfies 
$A(x) = \sum_{k=0}^{\infty} k! \big(x A(x)\big)^{k}$ and hence is given by 
\htmladdnormallink{A088368}{http://www.research.att.com:80/cgi-bin/access.cgi/as/njas/sequences/eisA.cgi?Anum=A088368},
and we deduce a moderately efficient recurrence relation. In \S4 we 
define the noncrossing partition transform on integer sequences and give some 
examples.
In \S 5 we adapt the method of \S3 to obtain an analogous recurrence for 
the ``list of noncrossing lists'' case.

\vspace*{5mm}

\section{Lists Of Noncrossing Sets 
$\longleftrightarrow$ Sets Of Lists}
It is easy to count partitions of $[n]$ into 
sets of $k$ lists: start with all $n!$ permutations of 
$[n]$; then for each one 
choose $k-1$ of the $n-1$ spaces between its entries to split it into a 
list of $k$ nonempty lists. This yields all partitions of $[n]$ into 
lists of $k$ lists and shows that there are $n!\binom{n-1}{k-1}$ of them. 
Finally, to count partitions $[n]$ into sets (rather than lists) of 
$k$ 
lists, divide by $k!$. The result is $\frac{n!}{k!}\binom{n-1}{k-1}$, the so-called Lah 
number $L(n,k)$
(\htmladdnormallink{A105278}{http://www.research.att.com:80/cgi-bin/access.cgi/as/njas/sequences/eisA.cgi?Anum=A105278}). 

To count the lists of noncrossing sets,
first recall the well known bijection (essentially due to Prodinger 
\cite{prod83}) from Dyck $n$-paths to noncrossing partitions of $[n]$ 
illustrated below.
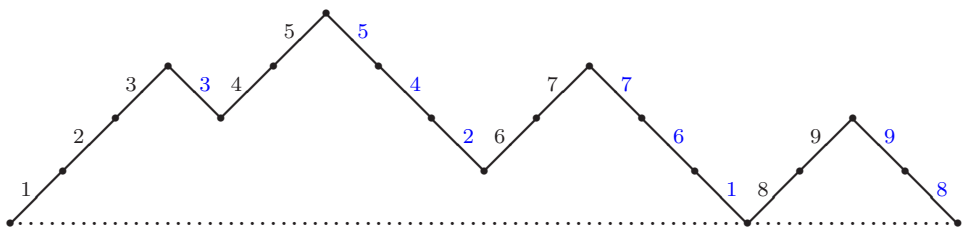
\Einheit=0.7cm
\[
\Pfad(-9,0),333433444334443344\endPfad
\SPfad(-9,0),111111111111111111\endSPfad
\DuennPunkt(-9,0)
\DuennPunkt(-8,1)
\DuennPunkt(-7,2)
\DuennPunkt(-6,3)
\DuennPunkt(-5,2)
\DuennPunkt(-4,3)
\DuennPunkt(-3,4)
\DuennPunkt(-2,3)
\DuennPunkt(-1,2)
\DuennPunkt(0,1)
\DuennPunkt(1,2)
\DuennPunkt(2,3)
\DuennPunkt(3,2)
\DuennPunkt(4,1)
\DuennPunkt(5,0)
\DuennPunkt(6,1)
\DuennPunkt(7,2)
\DuennPunkt(8,1)
\DuennPunkt(9,0)
\Label\o{ \textrm{{\scriptsize 1}}}(-8.7,0.3)
\Label\o{ \textrm{{\scriptsize 2}}}(-7.7,1.3)
\Label\o{ \textrm{{\scriptsize 3}}}(-6.7,2.3)
\Label\o{ \textrm{{\scriptsize 4}}}(-4.7,2.3)
\Label\o{ \textrm{{\scriptsize 5}}}(-3.7,3.3)
\Label\o{ \textrm{{\scriptsize 6}}}(0.3,1.3)
\Label\o{ \textrm{{\scriptsize 7}}}(1.3,2.3)
\Label\o{ \textrm{{\scriptsize 8}}}(5.3,0.3)
\Label\o{ \textrm{{\scriptsize 9}}}(6.3,1.3)
\blue{
\Label\o{ \textrm{{\scriptsize 3}}}(-5.3,2.3)
\Label\o{ \textrm{{\scriptsize 5}}}(-2.3,3.3)
\Label\o{ \textrm{{\scriptsize 4}}}(-1.3,2.3)
\Label\o{ \textrm{{\scriptsize 2}}}(-0.3,1.3)
\Label\o{ \textrm{{\scriptsize 7}}}(2.7,2.3)
\Label\o{ \textrm{{\scriptsize 6}}}(3.7,1.3)
\Label\o{ \textrm{{\scriptsize 1}}}(4.7,0.3)
\Label\o{ \textrm{{\scriptsize 9}}}(7.7,1.3)
\Label\o{ \textrm{{\scriptsize 8}}}(8.7,0.3) }
\Label\u{ \textrm{\small number the upsteps left to right,}}(0,-0.5)
\Label\u{ \textrm{\small label each downstep with the number on its 
matching 
upstep,}}(0,-1.3)
\Label\u{ \textrm{\small form the partition of $[n]$ whose blocks 
are the labels on the descents.}}(0,-2.1)
\]

\vspace*{3mm}

\noindent The Dyck path shown thus corresponds to the noncrossing partition 
3-542-761-98 (in a standard form: entries decreasing in each block and 
blocks listed in increasing order of their first entries).
This bijection sends \# peaks in the Dyck path to \# blocks in the 
partition. Since a noncrossing partition of $[n]$ with $k$ blocks 
gives rise to $k!$ lists of sets, partitions of $[n]$ into 
lists of noncrossing sets correspond to peak-labeled Dyck $n$-paths 
where peak-labeled means the peaks are labeled $1,2,3,\ldots$ in some 
order. Now the Narayana 
number $N(n,k) = \frac{1}{n}\binom{n}{k}\binom{n}{k-1}$ 
(\htmladdnormallink{A001263}{http://www.research.att.com:80/cgi-bin/access.cgi/as/njas/sequences/eisA.cgi?Anum=A001263})
is known to count Dyck $n$-paths with $k$ 
peaks, and so the number of peak-labeled Dyck $n$-paths with $n+1-k$ 
peaks is $(n+1-k)!N(n,k)$, which simplifies to the Lah number 
$L(n,k)$ mentioned above.

Summing over $k$ in the two preceding paragraphs yields the equivalence of the section title. However, we wish to 
show this equivalence directly by
giving a  bijection from peak-labeled Dyck $n$-paths with $k$ 
peaks to partitions of $[n]$ into sets of $n+1-k$ lists. Using the 
Dyck path above as a working example (with $n=9$ and $k=4$), begin by 
\emph{prepending an upstep}. Record the peak labels, ascent lengths, 
and descent lengths in left to right order as shown on the left below. 

\noindent\begin{tabular}{ccccc}
labels & 3 & 1 & 4 & 2  \\ \hline
ascents & 4 & 2 & 2 & 2  \\
descents & 1 & 3 & 3 & 2   \\
\end{tabular}\ $\overset{\raisebox{1mm}{ \textrm{{\scriptsize 
rotate}} }}{ \underset{ \textrm{{\scriptsize 
columns}} }{\longrightarrow} }$ \ 
\begin{tabular}{cccc}
 2 & 3 & 1 & 4    \\ \hline
 2 & 4 & 2 & 2    \\
 2 & 1 & 3 & 3    \\
\end{tabular}\ $\overset{\raisebox{1mm}{ \textrm{{\scriptsize 
delete}} }}{ \underset{ \textrm{{\scriptsize 
last col}} }{\longrightarrow} }$ \ 
\begin{tabular}{ccc}
 2 & 3 & 1   \\ \hline
 2 & 4 & 2   \\
 2 & 1 & 3   \\
\end{tabular}\ $\overset{\raisebox{1mm}{ \textrm{{\scriptsize 
partial}} }}{ \underset{ \textrm{{\scriptsize 
sums}} }{\longrightarrow} }$  \
\begin{tabular}{ccc}
 2 & 3 & 1     \\ \hline
 2 & 6 & 8     \\
 2 & 3 & 6    \\
\end{tabular}

The arrows illustrate the following steps: (i) cyclically rotate the columns 
so that the largest peak label is last, (ii)
drop the last column, and (iii) form partial sums of the bottom two rows.
Now form $[n]\,\backslash\,\{\textrm{middle row}\} = 
[9]\,\backslash\,\{2,6,8\}=\{1,3,4,5,7,9\}$---these numbers will be the first 
entries of the lists---and $[n]\,\backslash\,\{\textrm{bottom row}\} = 
[9]\,\backslash\,\{2,3,6\}=\{1,4,5,7,8,9\}$ and apply the difference 
operator (leaving the first entry intact) to get $\{1,3,1,2,1,1\}$---these 
numbers will be the lengths of the lists. From their lengths and 
first entries, we now have 
partial lists
\Einheit=0.5cm
\[
\Label\o{1}(-6.5,0)
\Pfad(-7,0),1\endPfad
\Pfad(-7,1),1\endPfad
\Pfad(-7,0),2\endPfad
\Pfad(-6,0),2\endPfad
\Label\o{3}(-4.5,0)
\Pfad(-5,0),111\endPfad
\Pfad(-5,1),111\endPfad
\Pfad(-5,0),2\endPfad
\Pfad(-4,0),2\endPfad
\Pfad(-3,0),2\endPfad
\Pfad(-2,0),2\endPfad
\Label\o{4}(-0.5,0)
\Pfad(-1,0),1\endPfad
\Pfad(-1,1),1\endPfad
\Pfad(-1,0),2\endPfad
\Pfad(-0,0),2\endPfad
\Label\o{5}(1.5,0)
\Pfad(1,0),11\endPfad
\Pfad(1,1),11\endPfad
\Pfad(1,0),2\endPfad
\Pfad(2,0),2\endPfad
\Pfad(3,0),2\endPfad
\Label\o{7}(4.5,0)
\Pfad(4,0),1\endPfad
\Pfad(4,1),1\endPfad
\Pfad(4,0),2\endPfad
\Pfad(5,0),2\endPfad
\Label\o{9}(6.5,0)
\Pfad(6,0),1\endPfad
\Pfad(6,1),1\endPfad
\Pfad(6,0),2\endPfad
\Pfad(7,0),2\endPfad
\]
and all that remains is to fill in the blanks. This is done by 
arranging the missing numbers in the order of their associated labels 
as in the table following the last arrow above---thus 2,6,8 in the order 2,3,1 is 
6,8,2---and then inserting them left to right in the blank squares. The 
final result is
\Einheit=0.5cm
\[
\Label\o{1}(-6.5,0)
\Pfad(-7,0),1\endPfad
\Pfad(-7,1),1\endPfad
\Pfad(-7,0),2\endPfad
\Pfad(-6,0),2\endPfad
\Label\o{3}(-4.5,0)
\Label\o{6}(-3.5,0)
\Label\o{8}(-2.5,0)
\Pfad(-5,0),111\endPfad
\Pfad(-5,1),111\endPfad
\Pfad(-5,0),2\endPfad
\Pfad(-4,0),2\endPfad
\Pfad(-3,0),2\endPfad
\Pfad(-2,0),2\endPfad
\Label\o{4}(-0.5,0)
\Pfad(-1,0),1\endPfad
\Pfad(-1,1),1\endPfad
\Pfad(-1,0),2\endPfad
\Pfad(-0,0),2\endPfad
\Label\o{5}(1.5,0)
\Label\o{2}(2.5,0)
\Pfad(1,0),11\endPfad
\Pfad(1,1),11\endPfad
\Pfad(1,0),2\endPfad
\Pfad(2,0),2\endPfad
\Pfad(3,0),2\endPfad
\Label\o{7}(4.5,0)
\Pfad(4,0),1\endPfad
\Pfad(4,1),1\endPfad
\Pfad(4,0),2\endPfad
\Pfad(5,0),2\endPfad
\Label\o{9}(6.5,0)
\Pfad(6,0),1\endPfad
\Pfad(6,1),1\endPfad
\Pfad(6,0),2\endPfad
\Pfad(7,0),2\endPfad
\]
giving the lists in increasing order of their first entries. We leave 
the interested reader to verify that the mapping is invertible; an appeal 
to the cycle lemma (see e.g., \cite[pp.\,359--360]{gkp}) will be needed to determine the appropriate cyclic 
rotation.

\vspace*{5mm}

\section{Sets Of Noncrossing Lists} Sequence 
\htmladdnormallink{A088368}{http://www.research.att.com:80/cgi-bin/access.cgi/as/njas/sequences/eisA.cgi?Anum=A088368}
is defined by the generating function equation $A(x) = \sum_{k=0}^{\infty} 
k! \big(x A(x)\big)^{k}$. We will show that this sequence counts partitions 
of $[n]$ into sets of noncrossing lists. Let $\u(n)$ denote this 
set of partitions and $\u(n,k)$ the subset for which $n$ occurs in a 
list of length $k$. Set $u(n)=\v \u(n) \v$ and $u(n,k)=\v \u(n,k) \v$; 
thus $u(n)=\sum_{k=1}^{n}u(n,k)$. For a partition in $\u(n,k)$ the 
entries in the list containing $n$ split $[n]$ into a sequence of subintervals 
$I_{1},I_{2},\ldots,I_{k}$ of lengths, say, 
$a_{1},a_{2},\ldots,a_{k}\ (a_{i}\ge 1,\ 1\le i \le k,\ 
\sum_{i=1}^{k}a_{i}=n$). Thus with $n=8$ the list 3,8,4 yields 
$[1,2,3],\:[4],\:[5,6,7,8]$. Set $J_{i}=I_{i} \backslash \{a_{i}\},\ 
1\le i \le k$.  The remaining lists are formed from entries of the 
$J_{i}$'s and since no crossovers are allowed between these lists (the 
noncrossing property would be violated), we are restricted to 
partitioning each $J_{i}$ into a set of noncrossing lists. This can be 
done in $u(b_{1})u(b_{2})\cdots u(b_{k})$ ways where $b_{i}=a_{i}-1,\ 
1\le i \le k$. Clearly, $b_{i}\ge 0$ and $\sum_{i=1}^{k}b_{i}=n-k$. Thus, 
\begin{equation}
u(n,k)=k! \sum_{(b_{1},\ldots,b_{k})}u(b_{1})u(b_{2})\cdots u(b_{k})    
    \label{eq:1}
\end{equation}
where the sum is taken over all nonnegative $k$-tuples 
$(b_{1},b_{2},\ldots,b_{k})$ whose sum is 
$n-k$ (weak compositions of $n-k$ into $k$ parts), 
and the $k!$ factor serves to order the block containing $n$.
Let $U(x)=\sum_{n\ge 0}u(n)x^{n}$ with $u(0):=1$.
Then the right hand side in (\ref{eq:1}) is $ [x^{n-k}]k!\big(U(x)\big)^{k}
=[x^{n}]k!\big(xU(x)\big)^{k}$. Multiply 
by $x^{n}$ and sum over $n$ and $k$ to get $U(x) = \sum_{k=0}^{\infty} 
k! \big(x U(x)\big)^{k}$, as claimed. 

The recurrence (\ref{eq:1}) is not efficient: there are $\binom{n-1}{k-1}$ nonnegative 
$k$-tuples whose sum is $n-k$. Thus to compute $u(n)$ using  
 (\ref{eq:1}) involves a sum over 
$\sum_{k=1}^{n}\binom{n-1}{k-1}=2^{n-1}$ 
terms.  It is easy, however, to reduce it to a sum over (integer) partitions of 
$n$, a set whose size, turning the famous 
Hardy-Rademacher-Ramanujan formula 
into round figures, is approximately $\frac{1}{7n}13^{\sqrt{n}} <\!\!< 
\,2^{n-1}$. Count frequencies in a weak composition 
$(b_{1},b_{2},\ldots,b_{k})$ of $n-k$ indexing a summand in 
(\ref{eq:1}) to get a weak partition
$0^{p_{0}}\,1^{p_{1}}\,2^{p_{2}}\,\ldots\,(n-1)^{p_{n-1}}$ (in frequency-of-parts 
form) of $n-k$ into $k$ parts; thus $\sum_{i=0}^{n-1}p_{i}=k$  
 and $\sum_{i=0}^{n-1}ip_{i}=n-k$. 
Each such weak partition of $n-k$ comes from 
$\binom{p_{0}+p_{1}+\cdots+p_{n-1}}{p_{0},\,p_{1},\,\ldots,\,p_{n-1}}$ 
weak compositions  of $n-k$, 
all of which make the same contribution $k!\prod_{i=1}^{n-1}u(i)^{p_{i}} =
(p_{0}+p_{1}+\cdots+p_{n-1})!\prod_{i=1}^{n-1}u(i)^{p_{i}}$ to 
the sum $u(n,k)$. Furthermore, all such weak partitions (regardless of $k$) 
arise by subtracting 1 from each part of a partition of $n$. These 
observations translate into a faster recurrence for $u(n)$:
\begin{equation}
    \begin{array}{ccl}
    u(0) & = & 1,\textrm{\quad and for $n\ge 1$} \\[1ex]
   u(n) & = & \sum (p_{1}+\cdots+p_{n})!
\displaystyle{\binom{p_{1}+\cdots+p_{n}}{p_{1},\,\ldots,\,p_{n}}}u(0)^{p_{1}}u(1)^{p_{2}}\cdots u(n-1)^{p_{n}}
\end{array}
\label{eq:2}
\end{equation}

\noindent where the sum is over all partitions 
$1^{p_{1}}\,2^{p_{2}}\,\ldots \,n^{p_{n}}$ of $n$.

\vspace*{5mm}

\section{The ``Noncrossing Partition'' Transform}
A closer look at the previous section suggests a transform on integer sequences 
and a combinatorial interpretation of it. For a sequence $(a_{k})_{k\ge 
0}$ with $a_{0}=1$, define $(b_{k})_{k\ge 0}$ by
\begin{equation}
    B(x)=\frac{1}{x}\left(\frac{x}{A(x)}\right)^{\langle -1 \rangle}
    \label{eq:tran1}
\end{equation}
where $A$ and $B$ are the ordinary generating functions for $(a_{k})$ 
and  $(b_{k})$ respectively, and ${}^{\langle -1 \rangle}$ denotes compositional 
inverse (reversion of series). Equivalently, $B(x)$ is the unique 
power series satisfying
\begin{equation}
    \sum_{k\ge 0}a_{k}\big(xB(x)\big)^{k} = B(x).
    \label{eq:tran2}
\end{equation}

The case $a_{k}=k!$ was treated in the previous section, where $b_{k}$ 
was then shown to count partitions of $[k]$ into noncrossing lists. 
The argument readily generalizes, however, from $k!$ to arbitrary 
$a_{k}$ (subject to $a_{0}=1$) to establish the following interpretation 
for $b_{k}$. If $a_{k}$ counts a class of combinatorial 
configurations, say $A$-structures, on $k$-sets, then $b_{k}$ counts 
the configurations obtained thusly:

\emph{Partition the set $[k]$ into noncrossing blocks and then put an 
$A$-structure on each block.}

For this reason, we call the transform $(a_{k}) \rightarrow (b_{k})$ 
defined by (\ref{eq:tran1}) the \emph{noncrossing partition} transform. 
Note that, necessarily, $b_{0}=1$ and $b_{1}=a_{1}$. Here are a few 
examples (in all cases, $a_{0}=1$ and $b_{0}=1$).

\vspace*{-6mm}

\[
\begin{array}{c|c}
    (a_{k})_{k\ge 1} &  (b_{k})_{k\ge 1} \\ \hline
    1 &   \textrm{\htmladdnormallink{$C_{k}$}{http://www.research.att.com:80/cgi-bin/access.cgi/as/njas/sequences/eisA.cgi?Anum=A000108}} \\ 
    2^{k} & \frac{2^{k}}{k+1}\binom{2k}{k} \\ 
    \frac{1}{k+1}\binom{2k}{k} & \frac{1}{2k+1}\binom{3k}{k} \\ 
    \frac{1}{2}\binom{2k}{k}  & 2^{k-1}C_{k} 
\end{array} \hspace*{10mm}
\begin{array}{c|c}
    (a_{k})_{k\ge 1} &  (b_{k})_{k\ge 1} \\ \hline
    \textrm{\htmladdnormallink{$F_{k-1}$}{http://www.research.att.com:80/cgi-bin/access.cgi/as/njas/sequences/eisA.cgi?Anum=A000045}}
    &  \textrm{\htmladdnormallink{$\triangle$-free dissections }{http://www.research.att.com:80/cgi-bin/access.cgi/as/njas/sequences/eisA.cgi?Anum=A046736}} \\ 
    2^{k-1} & 
    \textrm{\htmladdnormallink{little Schr\"{o}der \# }{http://www.research.att.com:80/cgi-bin/access.cgi/as/njas/sequences/eisA.cgi?Anum=A001003}} \\
    \textrm{\htmladdnormallink{little Schr\"{o}der \#}{http://www.research.att.com:80/cgi-bin/access.cgi/as/njas/sequences/eisA.cgi?Anum=A001003}}
    &  \textrm{\htmladdnormallink{``blobs''}{http://www.research.att.com:80/cgi-bin/access.cgi/as/njas/sequences/eisA.cgi?Anum=A003168}} \\
   C_{k-1} &  \textrm{\htmladdnormallink{ big Schr\"{o}der \#}{http://www.research.att.com:80/cgi-bin/access.cgi/as/njas/sequences/eisA.cgi?Anum=A006318}}
\end{array}
\]

\vspace*{1mm}

\centerline{{\small The noncrossing partition transform $(a_{k}) \rightarrow (b_{k})$}} 

\vspace*{2mm}

\noindent In particular, if $(a_{k})$ counts permutations of $[k]$ with some 
property, then $b_{k}$ counts partitions of $[k]$ into noncrossing 
lists, each of which has the property in question. For example, since 
321-avoiding permutations are counted by the Catalan numbers, we see 
from the table above that the number of partitions of $[k]$ into noncrossing 
321-avoiding permutations is $\frac{1}{2k+1}\binom{3k}{k}$.

\vspace*{5mm}

\section{Lists Of Noncrossing Lists}
It is possible to use the decomposition of the block containing $n$ 
as in the Section 3 to obtain recurrence relations for the 
number of partitions of $[n]$ into lists of noncrossing lists. (We 
continue to use the descriptive term ``block'' but now it means a 
list rather than a set.) Here, however, 
the factor $u(i)$ in the product on the right side of (\ref{eq:2}) will 
count \emph{lists} of blocks, and so it will be necessary to remove the order 
on each such list of blocks, throw the block containing $n$ into the 
mix, and then re-order the whole lot. This requires keeping track of 
the number of blocks. So let $u(n,j)$ denote the number of partitions of 
$[n]$ into a list of $j$ noncrossing lists.

In a partition $\mbf{\Pi}$ of $[n]$ into lists of noncrossing lists, 
let $k$ denote the length of the block containing $n$. As in 
\S3, this block induces a decomposition of $[n]$ into intervals  
$I_{1},I_{2},\ldots,I_{k}$ whose terminal points form the block. 
Deleting the endpoints gives a list of intervals 
$J_{1},J_{2},\ldots,J_{k}$ each of which is a union of blocks in  
$\mbf{\Pi}$. Let $(a_{i})_{i=1}^{k}$ denote the lengths of the 
$I_{i}$ taken in decreasing order so that $\pi=(a_{i})_{i=1}^{k}$ is 
an integer partition of $n$. Set $b_{i}=a_{i}-1$ and suppose the 
first $r$ $b_{i}$'s are positive. Then for $i=1,2,\ldots,r$ the 
$J$ interval corresponding to $b_{i}$ is the union of some number 
$c_{i}\ (1\le c_{i} \le b_{i})$ of blocks in the original partition 
$\mbf{\Pi}$. The total number of blocks in $\mbf{\Pi}$ is then 
$1+\sum_{i=1}^{r}c_{i}$. Write the integer partition $\pi$ of $n$ in 
frequency-of-parts form $1^{p_{1}}\,2^{p_{2}}\,\ldots\,(n)^{p_{n}}$; 
thus $\sum ip_{i}=n$ and $\sum p_{i}=k$.

Thus summands in the recursive sum for $u(n,j)$ are indexed by 
configurations of the form
\begin{equation}
\begin{array}{ccccccc}
    b_{1} & b_{2} & \ldots & b_{r} & 0 & \ldots & 0  \\
    c_{1} & c_{2} & \ldots & c_{r} &  &  & 
\end{array}
   \label{eq:3}
\end{equation}
where  $1\le c_{i}\le b_{i}$ for $1\le i \le r$, and the top row is 
$\pi-1$ (entrywise) for some integer partition $\pi=(a_{i})_{i=1}^{k}$ of 
$n$, and $1+\sum_{i=1}^{r}c_{i}=j$.

Then, with the sum taken over these configurations,
\[
\begin{array}{lcl}
u(n,j)=\sum & \quad  k! & \textrm{{\small [\,permute entries of block 
containing 
$n$\,] $\times$}} \\
 &  \binom{p_{1}+\ldots+p_{n}}{p_{1},\ldots,p_{n}} & \textrm{{\small 
 [\,permute the lengths
  $|I_{1}|,\ldots,|I_{k}|$\,] $\times$}} \\
  & \quad  j! & \textrm{{\small  [\,permute blocks in $\mbf{\Pi}$\,] $\times$}} \\
 & \frac{u(b_{1},c_{1})}{c_{1}!}\frac{u(b_{2},c_{2})}{c_{2}!} \cdots 
\frac{u(b_{r},c_{r})}{c_{r}!} & \textrm{{\small [\,the denominators eliminate the 
inter-block}}\\
 & &\, \textrm{{\small \  order captured by $u$\,]}}
\end{array}
\]

Equivalently,
\begin{equation}
  u(n,j)=\sum (p_{1}+\cdots+p_{n})! \binom{p_{1}+\cdots+p_{n}}{p_{1},\,\ldots,\,p_{n}}
   \binom{1+c_{1}+\cdots+c_{r}}{1,c_{1},\,\ldots,\,c_{r}} u(b_{1},c_{1})\cdots 
   u(b_{r},c_{r})
   \label{eq:6}
\end{equation}
and $u(n):=\sum_{j=1}^{n}u(n,j)$ gives the number of partitions of 
$[n]$ into lists of noncrossing lists. 

The sequence $\big(u(n)\big)_{n\ge 1}$ 
begins $(1,\ 4,\ 24,\ 184,\ 1680,\ 17592,\ 206472,\ 2674752,\ \ldots )$.
The total number of terms $t(n)$ in the sum for $u(n)$ in (\ref{eq:6}) is $\sum 
b_{1}b_{2}\cdots b_{r}$ taken over all the configurations in 
(\ref{eq:3}) above. Thus $t(n)$ is the sum of products of the nonzero entries 
in  $\pi-1$ taken over all partitions $\pi$ of 
$n$. The generating function $\sum_{n\ge 0}t(n)x^{n}$ with $t(0):=1$ is 
given by
\begin{equation}
 \frac{1}{1-x}\ \prod_{k\ge 2}\frac{1}{1-(k-1)x^{k}}=1 + x + 2\,{x^2} + 4\,{x^3} + 8\,{x^4} + 14\,{x^5} + 27\,{x^6} + 45\,{x^7} + 82\,{x^8}+\ldots.
   \label{eq:7}
\end{equation}
(Cf. 
\htmladdnormallink{A006906}{http://www.research.att.com:80/cgi-bin/access.cgi/as/njas/sequences/eisA.cgi?Anum=A006906}
for the sum of products of the entries taken over all partitions of $n$.)

The number of terms in the sum for $u(n,j)$ can be somewhat 
further reduced by collecting equal $u(b_{i},c_{i})$ factors: if $b$ occurs 
$j$ times among the $b_{i}$, then collecting equal $c$'s reduces the 
contribution of $b$ to the number of terms from a factor of $b^{j}$ to 
a factor of $\binom{b+j-1}{j}$. The generating function for the total number 
of terms 
thereby changes from (\ref{eq:7})  to 
\[
\frac{1}{1-x}\ \prod_{k\ge 2}\frac{1}{(1-x^{k})^{k-1}}=1 + x + 2\,{x^2} + 4\,{x^3} + 8\,{x^4} + 14\,{x^5} + 26\,{x^6} + 44\,{x^7} + 77\,{x^8}+\ldots .
\]

\vspace*{5mm}
\textbf{\large{Acknowledgments}}\quad I thank the referee for a 
careful reading of the paper and several helpful suggestions.

\noindent 2000 {\it Mathematics Subject Classification}: 05A15.

\noindent \emph{Keywords: } Set partitions, lists, noncrossing, cycle lemma
\bigskip
\hrule
\bigskip

\noindent (Concerned with sequences
\seqnum{A000108},
\seqnum{A000110},
\seqnum{A000262},
\seqnum{A000670},
\seqnum{A001263},
\seqnum{A002866},
\seqnum{A006906},
\seqnum{A088368},
and
\seqnum{A105278}.)


\begin{thebibliography}{99}
\bibitem{prod83} Helmut Prodinger, A correspondence between ordered 
trees and noncrossing partitions,
\emph{Discrete Math.}  \textbf{46} (1983), Issue 2, 205--206.

\bibitem{gkp} Ronald L. Graham, Donald E. Knuth, Oren Patashnik, 
\emph{Concrete Mathematics} 
(2nd edition), Addison-Wesley, 1994.

\end{thebibliography}
\end{document}

xxE(&Ù`.©.©HMonacoZÜ&@ðÿÿÿÿ&Ç2ð"H6žÿüO —rO —rÃÏZÊxxE&oˆŽ>MPSR
ô&¢‰øíÿÿ,&¢ªÐcolors